\documentclass[11pt,a4paper]{article}
\usepackage{amsmath,amssymb,amsthm,enumerate,hyperref,mathrsfs,mathpazo}
\usepackage[all]{xy}
\usepackage[british]{babel}
\usepackage[latin1]{inputenc}\usepackage[T1]{fontenc}
\usepackage[mathcal]{euscript}
\usepackage[hmargin=2cm,vmargin={2cm}]{geometry}

\newtheorem{theorem}{Theorem}[section]
\newcommand{\bt}{\begin{theorem}}
\newcommand{\et}{\end{theorem}}

\newtheorem{prop}[theorem]{Proposition}
\newcommand{\bp}{\begin{prop}}
\newcommand{\ep}{\end{prop}}

\newtheorem{lemme}[theorem]{Lemma}
\newcommand{\bl}{\begin{lemme}}
\newcommand{\el}{\end{lemme}}

\theoremstyle{definition}

\newtheorem{dfn}[theorem]{{Definition}}
\newcommand{\bd}{\begin{dfn}}
\newcommand{\ed}{\end{dfn}}

\theoremstyle{remark}

\newtheorem{rmq}[theorem]{Remark}
\newcommand{\br}{\begin{rmq}}
\newcommand{\er}{\end{rmq}}

\newtheorem{notation}[theorem]{Notation}
\newcommand{\bn}{\begin{notation}}
\newcommand{\en}{\end{notation}}

\newcommand{\bpf}{\begin{proof}}
\newcommand{\epf}{\end{proof}}

\newcommand{\zz}{\mathbb{Z}}
\newcommand{\nn}{\mathbb{N}}

\newcommand{\td}{\mathbb{D}}
\newcommand{\ta}{\mathbb{A}}
\newcommand{\te}{\mathbb{E}}
\newcommand{\tabc}{\mathbb{T}}
\newcommand{\D}{\mathcal{D}}
\newcommand{\di}{\mathbf{dim}}
\newcommand{\Hom}{\operatorname{Hom}\nolimits}
\newcommand{\tr}{\operatorname{tr}\nolimits}

\newcommand{\Ext}{\operatorname{Ext}\nolimits}
\newcommand{\id}{\operatorname{id}\nolimits}
\newcommand{\ul}{\underline}
\newcommand{\ppq}{\leqslant}
\newcommand{\pgq}{\geqslant}
\newcommand{\ot}{\otimes}

\newcommand{\mx}[1]{\begin{pmatrix}#1\end{pmatrix}} 
\newcommand{\smx}[1]{\left(\begin{smallmatrix}#1\end{smallmatrix}\right)} 
\newcommand{\rep}[1]{\langle{#1}\rangle} 
\newcommand{\epv}[2]{\varepsilon_{#2}^{(#1)}}
\newcommand{\set}[1]{\left\{ #1 \right\}}

\title{Combinatorial classification of piecewise hereditary algebras}
\author{Marcelo Lanzilotta,   Maria Julia Redondo  and Rachel Taillefer}
\date{}

\begin{document}
\maketitle

\begin{flushleft}


\medskip

\end{flushleft}

\begin{abstract} We use the characteristic polynomial of the Coxeter matrix of an algebra to complete the combinatorial classification of piecewise hereditary algebras which Happel gave in terms of the trace of the Coxeter matrix. We also give a cohomological interpretation of the coefficients (other than the trace) of the characteristic polynomial  of the Coxeter matrix  of any finite dimensional algebra with finite global dimension. 
\end{abstract}

\section{Introduction}

The main aim of this note is to complete Happel's combinatorial classification of piecewise hereditary algebras. 

An algebra $A $ is said to be \textbf{piecewise hereditary of type ${\mathcal H}$} if its bounded derived category is triangle equivalent to the bounded derived category of a hereditary abelian $k$-category ${\mathcal H}$ for a field $k$ which we assume to be algebraically closed. Such algebras have been much studied, see for instance \cite{H3,Hinv,HR,HRS,Lache,Ladkani,LeMeur} among others.

Happel and Reiten \cite{HR} have shown that ${\mathcal H}$ has a
tilting object and therefore $\D^b({\mathcal H})$ is triangle
equivalent to $\D^b(\Lambda)$ where $\Lambda$ is a  finite dimensional
hereditary algebra (\emph{ie.\/} a path algebra $k\vec{\Delta}$ for a
quiver $\vec{\Delta}$) or a canonical algebra (see \cite{Hinv}). In the
first case we say that $A$ is of type $k\vec{\Delta}$ and in the
second case we say that $A$ is of canonical type.

Recall that a
\textbf{canonical algebra}  is a one-point extension (see Definition \ref{oneptext}) of the path algebra of the quiver
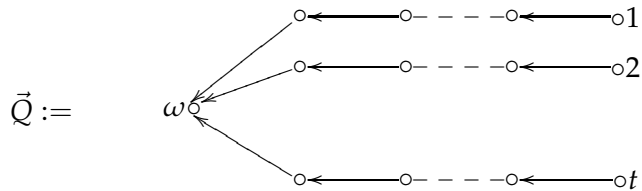
\[\xymatrix@R=0.2cm@C=1.2cm@M=0pt{&&\circ\ar[dddl]&\circ\ar[l]&\circ\ar@{--}[l]&\circ1\ar[l]\\&&&&&\\&&\circ\ar[dl]&\circ\ar[l]&\circ\ar@{--}[l]&\circ2\ar[l]\\ \vec{Q}:=&\omega\circ\\&&&&&\\&&&&&\\&&\circ\ar[uuul]&\circ\ar[l]&\circ\ar@{--}[l]&\circ t\ar[l]} \] 
by an indecomposable module $M$ given by the following representation
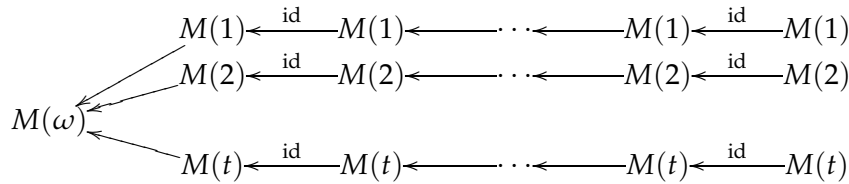
\[ \xymatrix@R=0.2cm@C=1.2cm@M=0pt{ 
&  M(1) \ar[ddl] & M(1) \ar[l]_{\id} & \cdots  \ar[l] & M(1) \ar[l] & M(1) \ar[l]_{\id} \\
&  M(2) \ar[dl]  & M(2) \ar[l]_{\id} & \cdots  \ar[l] & M(2) \ar[l] & M(2) \ar[l]_{\id} \\
M(\omega) & \\
&  M(t) \ar[ul]  & M(t) \ar[l]_{\id} & \cdots \ar[l] & M(t) \ar[l] & M(t) \ar[l]_{\id} } \]
such that $M(1), \dots, M(t)$ are pairwise different one-dimensional subspaces of $M(\omega)$, $\dim_k M(\omega)=2$ and $t \pgq 2$. 
 If we order the vertices from left to right and from top to bottom, the last vertex being $\omega$, then  $M$ has dimension vector $m:=\di_{k\vec{Q}}(M)=(1,\ldots,1,2)$

Happel has given the following combinatorial characterisation of some piecewise hereditary algebras, using the trace of the Coxeter matrix $\phi_A$ of $A$ (see Section \ref{prelim} for a definition), whose proof follows from \cite{H3} and  Happel's trace formula \cite{H}.

\bt[Happel]\label{Happelclassification} Let $A$ be a connected piecewise hereditary algebra over an algebraically closed field. Then 
\begin{enumerate}[(1)]
\item $A$ is of type $k \vec{\Delta}$ where the underlying graph $\Delta$ is not a tree if and only if  $\tr\phi_A>-1,$
\item $A$ is of canonical type with $t>3$ branches if and only if $\tr\phi_A<-1,$ 
\item $A$ is of canonical type  with $t=3$  branches or of type $k
  \vec{\Delta}$ where the underlying graph $\Delta$ is a tree if and only if $\tr\phi_A=-1$.
\end{enumerate}
\et

In this paper, we complete this classification, that is,  separate path algebra and canonical types when the trace of the Coxeter matrix is $-1$:

\bt\label{classification} Let $A$ be a connected  piecewise hereditary algebra over an algebraically
closed field with $n+1$ isomorphism classes of indecomposable projective modules and
$\tr\phi_A=-1$. Let $\chi_A(x)=\sum_{i=0}^{n+1}\lambda^A_ix^i$ be the characteristic polynomial of the Coxeter matrix $\phi_A$ of $A.$ Then $A$ is of canonical type (with $t=3$ branches) if
and only if one of the three following sets of conditions holds:
\begin{enumerate}[(i)]
\item $\lambda_{n-1}^A=0,$ $\lambda_{n-2}^A=-1=\lambda_{n-3}^A$ and $\lambda_{n-4}^A=0.$
\item $\lambda_{n-1}^A=0=\lambda_{n-2}^A$ and $\lambda_{n-\ell}^A\ppq -1$ for some $\ell\pgq 3.$
\item $\lambda_{n-1}^A=1$ and $\lambda_{n-\ell}^A\ppq0$ for some $\ell\pgq2.$
\end{enumerate}
\et

The proof of this theorem is the object of Section
\ref{sectionclassification}.

Finally, in Section \ref{cohomologicalinterpretation}, we give a cohomological interpretation of the coefficients of the characteristic polynomial of $\phi_A$ thus extending Happel's trace formula \cite{H}. 

\section{Preliminaries}\label{prelim}

Assume throughout that $k$ is an algebraically closed field and that all algebras are finite dimensional, basic, connected $k$-algebras and have finite global dimension.

\sloppy Let $A$ be an algebra and let $P(1),\ldots,P(n)$ be a set of
representatives of the isomorphism classes of indecomposable
projective left $A$-modules. For a finite dimensional left $A$-module $M$ we denote by $\di_A(M)\in\zz^n$ the dimension vector of $M$: its $i$th component is $\dim_k\Hom_A(P(i),M).$ Let $C_A$ denote the Cartan matrix of $A$, that is, the $n\times n$ matrix whose $j$th column is the transpose of $\di_AP(j).$ Since $A$ has finite global dimension, it is well known that $C_A$ is invertible over $\zz.$ We may therefore consider the following:
\begin{itemize}
\item The Coxeter matrix $\phi_A=-C_A^{-t}C_A$ of $A.$
\item\sloppy The Euler form associated to $A,$ $\rep{-,-}_A:\zz^n\times\zz^n\rightarrow\zz$ defined by $\rep{x,y}_A=xC_A^{-t}y.$ It is known that for two finite dimensional $A$-modules $X$ and $Y$ such that $X$ has finite projective dimension or $Y$ has finite injective dimension we have $\rep{\di_AX,\di_AY}_A=\sum_{i\pgq0}\dim_k\Ext^i_A(X,Y)$ (see \cite{R}). 
\end{itemize}

Our main object of study is the \textbf{Coxeter polynomial} of $A$, that is, the characteristic polynomial  $\chi_A(x)=\det(xI_n-\phi_A)$ of the Coxeter matrix of $A$ and in particular its coefficients. The first non-trivial coefficient, that is, the trace, was studied by Happel in \cite{H}, who proved that  $\tr\phi_A=-\rep{\di_{A^e}A,\di_{A^e}A}$ where $A^e=A\ot_k A^{op}$ is the enveloping algebra of $A.$

We are interested in this paper in some of the other coefficients of $\chi_A.$ 

\section{Classification of piecewise hereditary algebras}\label{sectionclassification}

The aim of this section is to prove Theorem \ref{classification}. We begin with a few
comments and results that are used in the proof.

\br\label{rkrestrictcases} Let $\mathscr{C}$ be a canonical algebra with $t$ branches. If $p_i-1>0$ is the number of vertices in the $i$th branch, then
  the sign of $\delta_p:=t-2-\sum_{i=1}^t\frac{1}{p_i}$ gives us
  information on the representation type of $\mathscr{C}$ (see \cite{GL}, in particular Remark 5.4 and Proposition 5.5): if $\delta_p<0$ then
  $\mathscr{C}$ is of domestic type and is derived equivalent to the module
  category of a tame hereditary algebra (in fact for $t=2$, in which case $\delta_p<0$, the algebra $\mathscr{C}$ is hereditary by definition), if $\delta_p=0$ then $\mathscr{C}$ is
  of tubular type (studied by Ringel in \cite{R}) and if $\delta_p>0$ then $\mathscr{C}$ is
  of wild type. 

 Therefore, in order to distinguish piecewise hereditary algebras of canonical type which are not of hereditary type, we need only consider the case $\delta_p \pgq 0$; necessarily we then have $t \pgq 3$. The case $t>3$ is characterised by the trace (Theorem \ref{Happelclassification}), so we may now assume that $t=3.$ In this case, $\sum_{i=1}^3\frac{1}{p_i}\ppq 1.$ 
\er

\bn When $t=3,$ we shall denote the quiver $\vec{Q}$ in the introduction by $\tabc_{a,b,c}$ where $a=p_1-1\ppq b=p_2-1\ppq c=p_3-1$ are the number of vertices on each branch in increasing order.
\en

The aim of this section is to  separate canonical type from  type $k \vec{\Delta}$ when the trace is equal to $-1.$ For this we shall use an inductive formula for the Coxeter polynomial, which gives the Coxeter polynomial of a one-point extension in terms of that of the original algebra. This was recently proved by Happel in \cite{Happel}. Recall:

\bd\label{oneptext} Let $B$ be a finite dimensional algebra and let $M$ be a
$B$-module. The \textbf{one-point extension} of $B$ by $M$ is the algebra
$B[M]=\mx{B&M\\0&k}$ (with usual matrix addition and
multiplication). 
\ed

The quiver of $B[M]$ contains the quiver of $B$ as a full subquiver
and there is an additional vertex, the \textbf{extension vertex}.

Note that when the global dimension of $B$ is finite, so is the global dimension of $B[M]$, so that no further restrictions on $B[M]$ are required.

Let $m$ denote the dimension vector of the $B$-module $M$,
\emph{i.e.\/}  $m=\di_BM.$  Then the
Cartan matrix of ${B[M]}$ is $\mx{1&0\\m^t&C_B}$ (if the extension
vertex has index $0$ so is placed ``before'' the others) and the
Coxeter matrix of $B[M]$ is $\mx{\rep{m,m}-1&-m\phi_{B}\\-C_B^{-t}m^t&\phi_B}$ where $\rep{-,-}=\rep{-,-}_B$ is the Euler form for $B$ to simplify notation.

Assume now that $A=B[M]$ is a one-point extension. 
Let $n$ be the number of isomorphism classes of indecomposable projectives for $B,$ so that the  number of isomorphism classes of indecomposable projectives for $A$ is $n+1.$ 
 Denote by $\chi_A(x)=\sum_{i=0}^{n+1}\lambda_i^Ax^i$ the  Coxeter polynomial of $A$ and by $\chi_B(x)=\sum_{i=0}^n\lambda_i^Bx^i$   the  Coxeter polynomial of $B.$

\bt\cite[Theorem 2.1]{Happel}\label{coefoneptext} With the notation
above, for any integer $\ell$ with $0\ppq \ell\ppq n$, we have \[\lambda_{n+1-\ell}^A=\lambda_{n-\ell}^B-(\rep{m,m}-1)\lambda_{n-(\ell-1)}^B-\sum_{i=1}^{\ell-1}\lambda_{n-\ell+i+1}^B\rep{m\phi_B^i,m}.\]
\et

The rest of the section is devoted to the proof of Theorem \ref{classification} and relies on several lemmas which we give now.

 Note that since the Coxeter polynomial is a derived invariant (see \emph{eg.\/} \cite{Happel}),  we need only consider the Coxeter polynomials of path algebras of trees and of canonical algebras with three branches.
Moreover, let $\vec{\Delta}_1$ and $\vec{\Delta}_2$ be two quivers with the same underlying tree graph $\Delta$. Let $A_1=k\vec{\Delta}_1$ and $A_2=k\vec{\Delta}_2$ be the corresponding path algebras. It is known by \emph{eg.\/} \cite{Hbook} and \cite{H2} that $A_1$ and $A_2$ are derived equivalent. Therefore the Coxeter polynomials of $A_1$ and $A_2$ are equal. We shall also denote them by $\chi_\Delta.$ 

The Coxeter polynomial of $k\vec{\ta}_{n+1}$ for the Dynkin graph $\ta_{n+1}$ is well known, it is equal to $\sum_{i=0}^{n+1}x^i$ (see \emph{eg.\/} \cite{Happel,S,B}).

\bl\label{stars} Consider the path algebra $k\vec{\tabc}_{a,b,c}$ with $a\ppq b\ppq c$. Then 
\[ \lambda^{k\vec{\tabc}_{a,b,c}}_{(a+b+c)-\ell}=\frac{(1-\ell)(2+\ell)}{2}\quad\mbox{ for $0\ppq \ell\ppq a.$} \]
Moreover, if $a=1$ then $\chi_{\tabc_{1,b,c}}(x)=x^{b+c+2}+\sum_{j=c+1}^{b+c+1}(j-b-c)x^j+\sum_{j=b+1}^c(1-b)x^j+\sum_{j=1}^b(2-j)x^j+1.$ In particular, for the Dynkin graph $\td_{n+1}=\tabc_{1,1,n-2}$ we get $\chi_{\td_{n+1}}(x)=x^{n+1}+x^n+x+1.$
\el

\br Note that the Coxeter polynomial of $\td_{n+1}$ has been computed elsewhere, see for instance  \cite{Happel,S}. Moreover, this lemma contains the cases of the Dynkin graphs $\te$ and of the Euclidean graphs $\widetilde{\te}$ which can also be found for instance in \cite{Happel,S}.
\er

\bpf
To compute the Coxeter polynomial of $k\vec{\tabc}_{a,b,c}$, we use a result of Boldt \cite[Corollary 3.2]{B} (see also \cite{S}) which gives\[ \chi_{\tabc_{a,b,c}}(x)= \left(\sum_{i=0}^ax^i\right)\left(\sum_{i=0}^{b+c+1}x^i\right)-x\left(\sum_{i=0}^{a-1}x^i\right)\left(\sum_{i=0}^bx^i\right)\left(\sum_{i=0}^cx^i\right)\] (since the Coxeter polynomials for the $\ta_\ell$ are known). It is easy to check that \[\left(\sum_{i=0}^ax^i\right)\left(\sum_{i=0}^{b+c+1}x^i\right)=\sum_{p=0}^a(p+1)x^p+\sum_{p=a+1}^{b+c+1}(a+1)x^p+\sum_{p=b+c+2}^{a+b+c+1}(a+b+c+1-p+1)x^p.\] Similarly, \[x\left(\sum_{i=0}^{a-1}x^i\right)\left(\sum_{i=0}^bx^i\right)\left(\sum_{i=0}^cx^i\right)=\left(\sum_{i=1}^{a}x^i\right)\left(\sum_{q=0}^{b-1}(q+1)x^q+\sum_{q=b}^{c}(b+1)x^q+\sum_{q=c+1}^{b+c}(b+c-q+1)x^q\right)\] in which the coefficient of $x^{a+b+c-\ell}$  is 
\[
\begin{cases}
\sum_{q=b+c-\ell}^{b+c}(b+c-q+1)=\frac{(\ell+1)(\ell+2)}{2}&\mbox{ if $0\ppq\ell\ppq a-1$}\\
\sum_{q=b+c-a}^{b+c-1}(b+c-q+1)=\frac{a(a+3)}{2}&\mbox{ if $\ell=a$}.\\
\end{cases}
\]

 Finally, the coefficient of $x^{p}$ in $\chi_{\tabc_{a,b,c}}(x)$ is $-\frac{(a+b+c-1-p)(a+b+c+2-p)}{2}$ if $b+c\ppq p\ppq a+b+c$ as required.
\epf

\bl\label{othertrees} Let $\Delta$ be a tree with $n+1$ vertices which is neither $\ta_{n+1}$ nor $\tabc_{a,b,c}$. Then $\lambda_{n-1}^\Delta\ppq -1.$
\el

\br This lemma contains the Euclidean cases $\widetilde{\td}_n$ whose Coxeter polynomials are known entirely (see \cite{Happel,S}).
\er

\bpf $\Delta$ is characterised by the fact that it has either two or more vertices of valency $3$ (the valency of a vertex is the number of edges connected to it) or at least one vertex of valency at least $4.$ Applying \cite[Theorem 4.8]{Happel} gives $\lambda_{n-1}^\Delta\ppq -1$ in all these cases.
\epf

Let $A=\mathscr{C}_{a,b,c}$ denote the canonical algebra which is a one-point extension $B[M]$ of the tree $B=k\tabc_{a,b,c}$ (with $a\ppq b\ppq c$) as defined in the introduction. Recall that we need only consider the cases where $\frac{1}{a+1}+\frac{1}{b+1}+\frac{1}{c+1}\ppq 1$. We have $m:=\di_B(M)=(1,\ldots,1,2)$.

We shall use Theorem \ref{coefoneptext} to compute coefficients of the Coxeter polynomial of $A$. We always have $\lambda_{n+1}^A=1$, $\lambda_{n}^A=-tr\phi_A=1$.

We first need the Coxeter matrix of $B;$ this can be determined using \cite[Proposition 3.1]{B}. Before we give it, we introduce some notation; for positive integers $p,q$ we set: 
\begin{itemize}
\item $J_p=\smx{0&0\\I_{p-1}&0}$ (a $p\times p$ matrix; $I_{p-1}$ denotes the identity matrix)
\item $K_{p,q}=\smx{1&\cdots&1\\0&\cdots&0\\\vdots&\cdots&\vdots\\0&\cdots&0}$ (a $p\times q$ matrix)
\item $\epv{q}{p}$ the $p$th vector in the canonical basis of $k^q$ ($p\ppq q$)
\item $v_q=\sum_{p=1}^q\epv{q}{p}=(1,\ldots,1).$
\end{itemize}
They satisfy the following rules: for any positive integers $p,q,r$ we have
\begin{alignat*}{2}
&v_qJ_q=v_q-\epv{q}{q}&\qquad&v_pK_{p,q}=v_q\\
&\epv{q}{p}J_q=
\begin{cases}
\epv{q}{p-1}\mbox{ if $p\pgq 2$}\\0\mbox{ if $p=1$}
\end{cases}&&\epv{q}{p}K_{q,r}=
\begin{cases}
0\mbox{ if $p\pgq 2$}\\v_r\mbox{ if $p=1$}.
\end{cases}
\end{alignat*}

With this notation, we have $m=(v_a,v_b,v_c,2v_1)$. Using  \cite[Proposition 3.1]{B} we have 
\[ \phi_B=\mx{J_a&K_{a,b}&K_{a,c}&K_{a,1}\\K_{b,a}&J_b&K_{b,c}&K_{b,1}\\K_{c,a}&K_{c,b}&J_c&K_{c,1}\\-K_{1,a}&-K_{1,b}&-K_{1,c}&-1}. \]

Moreover, for all $i\pgq 0$ we have $\rep{m\phi_B^i,m}=\rep{m\phi_B^i,m}^t=mC_B^{-1}(m\phi_B^i)^t$ and Boldt explains in \cite{B} how to compute $C_B^{-1}$; we get 
\[ mC_B^{-1}=(-\epv{a}{1},-\epv{b}{1},-\epv{c}{1},2\epv{1}{1}). \]We always have $\rep{m,m}=1$.

\bl\label{canonical12c} Let $A$ be the canonical algebra $\mathscr{C}_{1,2,c}$ with $c\pgq 5.$ Then $\lambda_{n-1}^A=0,$ $\lambda_{n-2}^A=-1=\lambda_{n-3}^A$ and $\lambda_{n-4}^A=0.$
\el

\bpf We have $m\phi_B=(0,\epv{2}{1},v_c-\epv{c}{c},v_1)$, $m\phi_B^2=(v_1,0,v_c-\epv{c}{c}-\epv{c}{c-1},v_1)$, $m\phi_B^3=(0,v_2,v_c-\epv{c}{c}-\epv{c}{c-1}-\epv{c}{c-2},v_1)$ and $m\phi_B^4=(v_1,v_2-\epv{2}{2},v_c-\epv{c}{c}-\epv{c}{c-1}-\epv{c}{c-2}-\epv{c}{c-3},v_1)$, so that $\rep{m\phi_B^i,m}=0$ for $1\ppq i\ppq 3$ and $\rep{m\phi_B^4,m}=-1.$ We know the coefficients of $\chi_B(x)$ from Lemma \ref{stars}: $\lambda_n^B=1,$ $\lambda_{n-1}^B=1,$ $\lambda_{n-2}^B=0,$ $\lambda_{n-3}^B=-1$, $\lambda_{n-4}^B=-1$ and $\lambda_{n-5}^B=-1.$ Therefore the formula in Theorem \ref{coefoneptext} gives the result.
\epf

\bl\label{canonical1bc} Let $A$ be the canonical algebra $\mathscr{C}_{1,b,c}$ with $b\pgq 3.$ Then $\lambda_{n-r}^A=0$ for $1\ppq r\ppq b-1$ and $\lambda_{n-b}^A<0.$
\el

\bpf We have $m\phi_B=(0,v_b-\epv{b}{b},v_c-\epv{c}{c},v_1)$ and for $2\ppq r\ppq b$, by induction,
\[ m\phi_B^r=(F_{r-2}v_1,F_{r-1}v_b-\sum_{k=0}^{r-2}F_{r-2-k}\epv{b}{b-k}-\epv{b}{b-r+1},F_{r-1}v_c-\sum_{k=0}^{r-2}F_{r-2-k}\epv{c}{c-k}-\epv{c}{c-r+1},F_{r-1}v_1) \]
where $F_r$ is the $r$th term in the Fibonacci sequence ($F_0=1=F_1$ and $F_{r+2}=F_{r+1}+F_r$ for $r\pgq 0$). Therefore $\rep{m\phi_B,m}=0,$ $\rep{m\phi_B^r,m}=-F_{r-2}$ for $2\ppq r\ppq b-1$ and $\rep{m\phi_B^b,m}=-F_{b-2}+1+\delta_{bc}.$ We also know by Lemma \ref{stars} that $\lambda_{n-i}^B=
\begin{cases}
1\mbox{ if $i=0$}\\ 2-i\mbox{ if $1\ppq i\ppq b+1$}
\end{cases}
$ so that by Theorem \ref{coefoneptext} we get 
\[ 
\begin{cases}
\lambda_{n-1}^A=0\\
\lambda_{n-r}^A=1-r+\sum_{i=2}^{r-1}(2+i-r)F_{i-2}+F_{r-2}\mbox{ if $2\ppq r\ppq b-1$}\\\lambda_{n-b}^A=-\delta_{bc}-b+\sum_{i=2}^{b-1}(2+i-b)F_{i-2}+F_{b-2}.
\end{cases}
 \] We use the well known formula $\sum_{i=p}^qF_i=F_{q+2}-F_{p+1}$ to get $\sum_{i=0}^qiF_i=\sum_{k=1}^q\sum_{i=k}^qF_i=qF_{q+2}-F_{q+3}+F_3.$ These finally give:
\[ 
\begin{cases}
\lambda_{n-r}^A=0\mbox{ if $1\ppq r\ppq b-1$}\\\lambda_{n-b}^A=-\delta_{bc}-1<0.
\end{cases}\qedhere
 \]
\epf

\bl\label{canonicalabc} Let $A$ be the canonical algebra $\mathscr{C}_{a,b,c}$ with $a\pgq 2.$ Then $\lambda_{n-r}^A=1$ for $1\ppq r\ppq a-1$ and $\lambda_{n-a}^A\ppq0.$
\el

\bpf For $1\ppq r\ppq a$, we prove by induction that the vector $m\phi_B^r$ is equal to  \[(2^{r-1}v_a-\sum_{j=0}^{r-2}2^{r-2-j}\epv{a}{a-j}-\epv{a}{a-r+1},2^{r-1}v_b-\sum_{j=0}^{r-2}2^{r-2-j}\epv{b}{b-j}-\epv{b}{b-r+1},2^{r-1}v_c-\sum_{j=0}^{r-2}2^{r-2-j}\epv{c}{c-j}-\epv{c}{c-r+1},2^{r-1}v_1)\] so that $\rep{m\phi_B^r,m}=-2^{r-1}$ if $1\ppq r\ppq a-1$ and $\rep{m\phi_B^a,m}=-2^{a-1}+1+\delta_{ab}+\delta_{ac}.$ Then using Theorem \ref{coefoneptext} and the coefficients of $\chi_B$ obtained previously in Lemma \ref{stars} we get 
\[ 
\begin{cases}
\lambda_{n-i}^A=1\mbox{ if $1\ppq i\ppq a-1$} \\
\lambda_{n-a}^A=-\delta_{ab}-\delta_{ac}\ppq0
\end{cases}
 \] (we use the relations $\sum_{\ell=1}^p2^{-\ell}(\ell+1)=3-2^{-p}(p+3)$ and $\sum_{\ell=1}^p\ell(\ell+1)2^{-\ell+1}=2^4-2^{1-p}(p^2+5p+8)$ obtained by differentiating  the identity $\sum_{\ell=1}^px^{\ell+1}=x^2\frac{x^p-1}{x-1}$ twice and evaluating at $x=2^{-1}$).
\epf

We now have all we need to finish the proof of the classification.

\bpf[Proof of Theorem \ref{classification}.] Consider the Coxeter
polynomial of $A=\mathscr{C}_{a,b,c}$ with $a\pgq 2$ (and
$n=a+b+c+1$). Then $\lambda_{n-1}^A=1$ and   using  Lemmas \ref{stars}
and \ref{othertrees} we see that the only tree that satisfies this is
$\ta_{n+1}.$ However, the coefficients of $x^{n-a}$ in the Coxeter
polynomials of  $\ta_{n+1}$ and $\mathscr{C}_{a,b,c}$  differ since it is $1$ for $\ta_{n+1}$ and nonpositive for $\mathscr{C}_{a,b,c}$ by Lemma \ref{canonicalabc}. Therefore the Coxeter polynomial of $\mathscr{C}_{a,b,c}$ is different from that of all trees.

Now consider the Coxeter polynomial of $\mathscr{C}_{1,2,c}$ with $c\pgq 2$ (and $n=c+4$).  Recall from Remark \ref{rkrestrictcases} that we need only consider the case where $1\pgq \frac{1}{a+1}+\frac{1}{b+1}+\frac{1}{c+1}=\frac{1}{2}+\frac{1}{3}+\frac{1}{c+1}$, \emph{ie.} $c\pgq5$, so that we assume $c\pgq 5.$ Using  Lemmas \ref{stars} and \ref{othertrees} we see that the only tree such that  the  coefficients of $x^{n+1-i}$ for $0\ppq i\ppq 4$ in its Coxeter polynomial  are the same as those for $\mathscr{C}_{1,2,c}$ is $\tabc_{1,2,c+1}$. But the coefficient of $x^{c}$ in the Coxeter polynomial of $\tabc_{1,2,c+1}$ is $-1$ whereas for $\mathscr{C}_{1,2,c}$ it is $0$   by Lemma \ref{canonical12c}. Therefore the Coxeter polynomial of $\mathscr{C}_{1,2,c}$ is different from that of all trees.

Finally consider  the Coxeter polynomial of $\mathscr{C}_{1,b,c}$ with $b\pgq 3$ (and $n=b+c+2$).   Using Lemmas \ref{stars} and \ref{othertrees} we see that the only tree   such that  the  coefficients of $x^{n+1-i}$ for $0\ppq i\ppq 3$ in its Coxeter polynomial are the same as those for $\mathscr{C}_{1,b,c}$ is $\td_{b+c+3}$. But the coefficient of $x^{c+2}$ is $0$ in the Coxeter polynomial of $\td_{b+c+3}$  and in that of $\mathscr{C}_{1,b,c}$ it is negative by Lemma \ref{canonical1bc}.  Therefore the Coxeter polynomial of $\mathscr{C}_{1,b,c}$ is different from that of all trees.
\epf

\section{Cohomological interpretation of the coefficients of the Coxeter polynomial}\label{cohomologicalinterpretation}

Happel proved  the following result:

\bt\cite{H}\label{trace} $\tr\phi_A=-\rep{\di_{A^e}A,\di_{A^e}A}.$
\et

We wish to do something similar for the other coefficients of the Coxeter polynomial. From now on, assume that the characteristic of $k$ is $0.$
We need an interpretation of these coefficients in terms of the entries of the Coxeter matrix:

\bp\label{charpolytrace} Let $\phi=(\phi_{ij})_{1\ppq i,j\ppq n}$ be a matrix, and let $\chi(x)=\det{(x\,\id -\phi)}$ be its characteristic polynomial. Write $\chi(x)=x^n+\lambda_{n-1}x^{n-1}+\cdots+\lambda_{1}x+\lambda_0.$
Then \[\lambda_{n-\ell}=\sum(-1)^{\sigma(\ul{p})}\alpha_{\ul{p}}\tr(\phi^{p_1})\cdots\tr(\phi^{p_r}),\] where the sum is taken over all partitions $\ul{p}=(p_1,\ldots,p_r)$ of $\ell$, $\sigma(\ul{p})=\sum_{i=1}^rp_i$ and  $\alpha_{\ul{p}} =\displaystyle{\frac{1}{p_1p_2\ldots p_r} \prod_{a=0}^\ell\frac{1}{\left(\#\set{i\,\mid\,p_i=a}\right)!}}$.
\ep

\bpf Without loss of generality, we may assume that the field $k$ is algebraically closed. Therefore the coefficient $\lambda_{n-\ell}$ of the characteristic polynomial of $\phi$ is the $\ell$th elementary symmetric polynomial $\sigma_\ell(\mu_1,\ldots,\mu_n)$ in the eigenvalues $\mu_1,\ldots,\mu_n$ of $\phi.$ This can be expressed in terms of the symmetric polynomials $S_k=\sum_{i=1}^n\mu_i^k=\tr(\phi^k)$ using Waring's formula, see for instance \cite[V.2]{LFA} or \cite[I.6]{MM}, which gives the expression above.
\epf
We now need to introduce some notation and results.

 Let $p_A(i)$ (\emph{resp.\/} $q_A(i)$, \emph{resp.\/} $e_A(i)$)
 denote the dimension vector of the $i$th indecomposable projective
 $A$-module $P(i)$ (\emph{resp.\/} the indecomposable injective $A$-module
 $Q(i)$, \emph{resp.\/} the $i$th simple $A$-module $S(i)$). Let
 $e_1,\ldots,e_n$ be the primitive orthogonal idempotents in $A$ such
 that $P(i)=Ae_i.$ Let $D$ denote the $k$-dual, \emph{ie.\/}
 $D=\Hom_k(-,k)$ and  let $\rep{-,-}_A$ denote the Euler form as before.

We have mentioned before that the transpose of $p_A(i)$ is the $i$th column of the Cartan matrix $C_A$. It is known that $q_A(i)$ is the $i$th row of $C_A$ (see \cite{R}).

We shall also need to work with bimodules. The indecomposable projective $A^e$-modules are the $Ae_i\ot e_j A.$ Set $e_{i,j}=e_i\ot e_j\in A^e,$ and let $S(i,j)$ be the corresponding simple module with $e_{A^e}(i,j)$ its dimension vector. 

We order the idempotents in the following way: \[e_{1,1},\ldots,e_{n,1},e_{1,2},\ldots,e_{n,2},\ldots,e_{1,n},\ldots,e_{n,n}.\]

 Happel recalled in \cite{H} the following results:
$\di_{A^e}A=(p_A(1),\ldots,p_A(n))$,  and $C_{A^e}^{-t}=C_A^{-1}\ot
C_A^{-t}$.

The dual $DA$ is an $A^e$-module. Its dimension vector is given by:
\begin{align*}
\dim_k\Hom_{A^e}(A^ee_{i,j},DA)&=\dim_k\Hom_{A^e}(A^ee_{i,j},\Hom_k(A,k))\\&=\dim_k\Hom_k(A\ot_{A^e}A^ee_{i,j},k)\\&=\dim_k D(A\ot_{A^e}A^ee_{i,j})=\dim_k D(Ae_{i,j})=\dim_k D(e_jAe_i)\\&=\dim_k(e_jAe_i)=\dim_k\Hom_A(P(j),P(i))=q_A(j)_i.
\end{align*} Therefore $\di_{A^e}DA=(q_A(1),\ldots,q_A(n))$.

We may now give a cohomological interpretation of the  trace of the powers
$\phi_A$ and hence of the coefficients of the Coxeter polynomial of $A$:

\bt Let $DA$ denote the dual of $A,$ viewed as a bimodule over $A.$ Then 
\begin{itemize}
\item $\tr (\phi_A^2)=\rep{\di_{A^e}DA,\di_{A^e}A}$.
\item If $k\pgq 3,$ then $(-1)^k\tr(\phi_A^k)$ is equal to \small
\[ \sum_{1\ppq v_1,\ldots,v_{k-1}\ppq n}\hspace*{-12pt}\rep{q_A(v_1),p_A(v_{k-1})}_A\rep{q_A(v_2),e_A(v_{1})}_A\ldots\rep{q_A(v_{k-2}),e_A(v_{k-3})}_A\rep{\di_{A^e}DA,e_{A^e}(v_{k-1},v_{k-2})}_{A^e} \]
\normalsize\end{itemize}
\et

\bpf We write $\rep{-,-}=\rep{-,-}_A$ to simplify notation. 
\begin{itemize}
\item We first prove that if $X\in\zz^n$ and $r\in\nn,$ $r\pgq 1,$ then $$(C_AC_A^{-t})^rX^t=\sum_{1\ppq u_1,\ldots,u_r\ppq n}\rep{q_A(u_1),X}\rep{q_A(u_2),e_A(u_1)}\ldots\rep{q_A(u_r),e_A(u_{r-1})}e_A(u_r)^t=:Y^t,$$ by induction on $r:$

If $r=1,$ we have $$C_AC_A^{-t}X^t=\mx{q_A(1)\\\vdots\\ q_A(n)}C_A^{-t}X^t=\mx{\rep{q_A(1),X}\\\vdots\\\rep{q_A(n),X}}=\sum_{u=1}^n\rep{q_A(u),X}e_A(u)^t.$$ 

Assume the result is true for $r;$ then 
\begin{align*}
(C_AC_A^{-t})^{r+1}X^t&=C_AC_A^{-t}Y^t=\sum_{u_{r+1}=1}^n\rep{q_A(u_{r+1}),Y}e_A(u_{r+1})^t\\
&=\sum_{1\ppq u_1,\ldots,u_r,u_{r+1}\ppq n}\rep{q_A(u_1),X}\rep{q_A(u_2),e_A(u_1)}\ldots\\&\hspace*{4cm}\rep{q_A(u_r),e_A(u_{r-1})}\rep{q_A(u_{r+1}),e_A(u_r)}e_A(u_{r+1})^t.
\end{align*}

\item  Set $C_A=(c_{ij})_{1\ppq i,j\ppq n}$ and $C_A^{-1}=(\zeta_{ij})_{1\ppq i,j\ppq n}.$ Then 
\begin{align*}
\tr{(\phi_A^2)}&=\sum_{1\ppq i,j\ppq n}\phi_{ij}\phi_{ji}
=\sum_{1\ppq i,j,r,s\ppq n}\zeta_{ri}c_{rj}\zeta_{sj}c_{si}=\sum_{i,r}\zeta_{ri}q_A(r)C_A^{-t}p_A(i)^t\\
&=\left(q_A(1),\ldots,q_A(n)\right) C_A^{-1}\ot C_A^{-t}\mx{p_A(1)^t\\\vdots\\p_A(n)^t}\\
&=\di_{A^e}(DA)\,C_{A^e}^{-t}(\di_{A^e}A)^t=\rep{\di_{A^e}DA,\di_{A^e}A}.
\end{align*}

\item Assume $k\pgq 3.$ Then 
\begin{align*}
\tr{(\phi_A^k)}&=\sum_{1\ppq j_1,\ldots,j_k\ppq n}\phi_{j_1j_2}\phi_{j_2j_3}\ldots\phi_{j_kj_1}= \sum_{\substack{
1\ppq j_1,\ldots,j_k\ppq n\\1\ppq i_1,\ldots,i_k\ppq n
}}
(-1)^k\zeta_{i_1j_1}c_{i_1j_2}\zeta_{i_2j_2}c_{i_2j_3}\ldots\zeta_{i_kj_k}c_{i_kj_1}\\
&= (-1)^k\sum_{1\ppq i_1,j_1\ppq n}\zeta_{i_1j_1}q_A(i_1)C_A^{-t}(C_AC_A^{-t})^{k-2}p_A(j_1)^t\\
&=(-1)^k\sum_{1\ppq v_1,\ldots,v_{k-1}\ppq n}\rep{q_A(v_1),p_A(v_{k-1})}\rep{q_A(v_{2}),e_A(v_{1})}\ldots \\&\hspace*{3.5cm}\rep{q_A(v_{k-2}),e_A(v_{k-3})}\di_{A^e}(DA)\, C_{A^e}^{-t}\, e_{A^e}(v_{k-1},v_{k-2})^t\\
&=(-1)^k\sum_{1\ppq v_1,\ldots,v_{k-1}\ppq n}\rep{q_A(v_1),p_A(v_{k-1})}\rep{q_A(v_2),e_A(v_{1})}\ldots\\&\hspace*{3.5cm}\rep{q_A(v_{k-2}),e_A(v_{k-3})}\rep{\di_{A^e}DA,e_{A^e}(v_{k-1},v_{k-2})}
\qedhere\end{align*}
\end{itemize}
\epf

\bigskip

\noindent\textsc{Marcelo Lanzilotta}\\
Centro de Matem\' atica (CMAT),  
Instituto de Matem\'atica y Estad\'\i stica Rafael Laguardia (IMERL), Universidad de la Rep\'ublica, Igu\'a 4225, C.P. 11400, Montevideo, Uruguay.\\\textit{E-mail address: marclan@cmat.edu.uy}\\\ \\
\textsc{Maria Julia Redondo}\\Instituto de Matem\'atica,
Universidad Nacional del Sur,
Av. Alem 1253, (8000) Bah\'\i a Blanca, Argentina.\\\textit{E-mail address:  mredondo@criba.edu.ar}\\\ \\
\textsc{Rachel Taillefer} (Corresponding author)\\Laboratoire LaMUSE, Universit\'e de Saint-Etienne, Facult\'e des Sciences, 23, rue du Dr. P. Michelon, 42023 Saint-Etienne Cedex 2, France.\\\textit{E-mail address: rachel.taillefer@univ-st-etienne.fr}

\end{document}